\documentclass[10pt,draft,twoside,leqno,a4paper]{article}
\usepackage{amsmath,amsthm,amsfonts,amssymb,verbatim}
\newtheorem{thm}{Theorem}[section]
\newtheorem{lemma}[thm]{Lemma}

\newtheorem{prop}[thm]{Proposition}

\numberwithin{equation}{section}

\newcommand{\R}{\mathbb{R}}
\newcommand{\N}{\mathbb{N}}
\newcommand{\Z}{\mathbb{Z}}
\newcommand{\C}{\mathbb{C}}
\newcommand{\I}{\mathbb{I}}
\newcommand{\Ed}{\mathcal E_\delta}
\newcommand{\Un}{U_N}
\newcommand{\e}{\mathrm{e}}
\newcommand{\Ln}{L_N}
\newcommand{\un}{u_n}
\newcommand{\unk}{u_{n_k}}
\newcommand{\ui}{u_\infty}

\newcommand{\eps}{\varepsilon}

\newcommand{\Ar}{A_R}
\newcommand{\disty}{d_y}
\newcommand{\distx}{d_x}
\newcommand{\Xd}{\hat X_\delta}
\newcommand{\Yd}{\hat Y_\delta}
\newcommand{\ro}{r_0}
\newcommand{\Bj}{\hat B_j}
\newcommand{\Bk}{\hat B_k}
\newcommand{\Qj}{\hat Q_j}
\newcommand{\phij}{\hat\varphi_j}
\newcommand{\psij}{\hat\psi_j}
\newcommand{\phik}{\hat\varphi_k}
\newcommand{\psik}{\hat\psi_k}
\newcommand{\gj}{\hat g_j}
\newcommand{\gk}{\hat g_k}
\newcommand{\hj}{\hat h_j}
\newcommand{\Cj}{\hat C_j}
\newcommand{\Ck}{\hat C_k}
\newcommand{\zd}{z_\delta}
\newcommand{\Fd}{F_\delta}
\newcommand{\Uj}{\hat U_j}
\newcommand{\Uk}{\hat U_k}
\newcommand{\Ld}{L_\delta}
\newcommand{\Lj}{\hat L_j}
\newcommand{\Lk}{\hat L_k}
\newcommand{\Lo}{\hat L_0}
\newcommand{\Sd}{S_\delta}
\newcommand{\Ap}{A'}
\newcommand{\phip}{\phi'}
\newcommand{\xp}{x'}
\newcommand{\ud}{u_\delta}
\newcommand{\hud}{\hat u_\delta}
\newcommand{\ek}{\hat{\underline e}_k}
\newcommand{\Rd}{R_\delta}
\begin{document}
\title{A shadowing lemma for
abelian Higgs vortices}
\author{Marta Macr\`{i}$^{(1)}$, Margherita Nolasco$^{(2)}$ and Tonia Ricciardi$^{(1)}$\\
\footnotesize{${}^{(1)}$Dipartimento di Matematica e Applicazioni,}
\footnotesize{Universit\`{a} di
Napoli Fe\-de\-ri\-co~II}\\ 
\footnotesize{Via Cintia, 80126 Napoli, Italy.}\\
\footnotesize{${}^{(2)}$Dipartimento di Matematica Pura ed Applicata,}
\footnotesize{Universit\`{a} di
L'Aquila}\\
\footnotesize{Via Vetoio, Coppito,}
\footnotesize{67010 L'Aquila, Italy.}\\
\small{\tt{macri@unina.it}}\\
\small{\tt{nolasco@univaq.it}}\\
\small{\tt{tonia.ricciardi@unina.it}}
}
\date{}
\maketitle
\begin{abstract}
We use a shadowing-type
lemma in order to analyze the singular, semilinear elliptic equation
describing static self-dual abelian Higgs
vortices.
Such an approach allows us to construct new solutions having an 
\textit{infinite} number of arbitrarily prescribed vortex points.
Furthermore, we obtain the precise asymptotic profile of the solutions
in the form of an approximate superposition rule, 
up to an error which is exponentially small.
\end{abstract}
\begin{description}
\item\small{ {\textsc{Key Words:} {Abelian Higgs model, elliptic equation, shadowing lemma}} }
\item \small{{\textsc{MSC 2000 Subject Classification:}} Primary 35J60; Secondary 58E15, 81T13}
\end{description}
\section{Introduction}
\label{sec:intro}
We consider the energy density for the static two-dimensional 
self-dual abelian Higgs model in the following form:
\begin{equation*}
\Ed(A,\phi)=\delta^2|\mathrm d A|^2+|D\phi|^2
+\frac{1}{4\delta^2}\left(|\phi|^2-1\right)^2,
\end{equation*}
where $A=A_1\mathrm{d}x_1+A_2\mathrm{d}x_2$, $A_1(x),A_2(x)\in\R$
is a gauge potential (a connection over a principal $U(1)$ bundle),
$\phi$, $\phi(x)\in\C$ is a Higgs matter field (a section over an
associated
complex line bundle), $D=\mathrm d-iA$
is the covariant derivative and $\delta>0$ is the coupling constant.
It corresponds to the two-dimensional Ginzburg-Landau energy density
in the so-called ``Bogomol'nyi limit", denoting the borderline between type I and type
II superconductors. 
In recent years, $\Ed$ has received considerable attention,
in view of both
its physical and geometrical interest, see, e.g.,  
\cite{GarciaPrada,HongJostStruwe,Stuart,Taubes,WangYang}
and the references therein.
\par
The smooth, finite action critical points for the action functional 
corresponding to $\Ed$ on $\R^2$ have been completely classified
by Taubes~\cite{JaffeTaubes,Taubes}. 
It is shown in \cite{JaffeTaubes} that such critical points
are completely determined by the distributional solutions to the elliptic problem
\begin{equation}
\label{eq:Taubes}
-\Delta u=\delta^{-2}(1-\e^u)-4\pi\sum_{j=1}^sm_j\delta_{p_j}
\qquad\mathrm{on\ }\R^2,
\end{equation}
which decay in the sense of the Sobolev space $H^1(\R^2)$ at infinity.
Here $s\in\N$, and for $j=1,2,\ldots,s$, $p_j\in\R^2$
are the vortex points, $m_j\in\N$ is the multiplicity of
$p_j$, $\delta_{p_j}$ is the Dirac measure at $p_j$.
By variational methods, Taubes proved that there exists a unique
solution
to \eqref{eq:Taubes} leading to a smooth, finite action critical point for 
the action functional of $\Ed$ on $\R^2$,
for any $s\in\N\cup\{0\}$, $p_j\in\R^2$ and $m_j\in\N$, 
$j=1,\ldots,s$,
and for any value of $\delta>0$.
Such a solution satisfies the topological constraint
$\int_{\R^2}F_{12}=2\pi\sum_{j=1}^sm_j$,
where $F_{12}=\partial_1A_2-\partial_2A_1$ is the magnetic field 
(the curvature of $A$).
\par
The case of infinitely many vortex points arranged on a
periodic lattice
has been considered in \cite{WangYang} 
and, in the more general setting of a compact Riemannian 2-manifold,
in \cite{GarciaPrada,HongJostStruwe}.
We say that the vortex points $p_j$, $j\in\N$ are doubly periodically
arranged in
$\R^2$ if there exists $s\in\N$ such that for any $k\in\N$, $k>s$
there exist
$j\in\{1,2\ldots s\}$ and $m,n\in\Z$ such that $p_k=p_j+m\underline
e_1+n\underline e_2$,
where $\underline e_1,\underline e_2$ are the  unit vectors in $\R^2$.
Similarly as in the previous case, denoting by $\Omega=\R^2/\Z^2$
the flat 2-torus, finite action critical points for the action of $\Ed$ on $\Omega$
correspond to distributional solutions to the problem
\begin{equation}
\label{eq:periodic}
-\Delta u=\delta^{-2}(1-\e^u)-4\pi\sum_{j=1}^sm_j\delta_{p_j}
\qquad\mathrm{on\ }\Omega,
\end{equation}
satisfying the topological constraint $\int_\Omega F_{12}=2\pi\sum_{j=1}^sm_j$.
It is shown in \cite{WangYang},
that
a unique solution for \eqref{eq:periodic} exists if
and only if $\delta\in(0,\pi^{-1})$. The asymptotics as $\delta\to0^+$
has been considered in \cite{HongJostStruwe,WangYang}.
\par
Our aim in this note is to show that
a shadowing lemma as introduced in the context of PDE's
by Angenent~\cite{Angenent},
see also \cite{Nolasco}, may be adapted in order to construct
solutions
to the following more general equation containing 
\textit{infinitely} many arbitrarily
prescribed vortex points:
\begin{equation}
\label{eq:infinite}
-\Delta u=\delta^{-2}(1-\e^u)-4\pi\sum_{j\in\N}m_j\delta_{p_j}
\qquad\mathrm{in\ }\R^2.
\end{equation}
Suitable modifications to the method described in \cite{Angenent}
are necessary, due to the singular sources appearing in \eqref{eq:infinite}.
We assume that the vortex points $p_j$,
$j\in\N$ are \textit{arbitrarily} distributed
in the plane, with the only constraint that 
\begin{equation}
\label{eq:constraints}
d:=\inf_{k\neq j}|p_j-p_k|>0
\qquad
\mathrm{and}\qquad m:=\sup_{j\in\N}m_j<+\infty. 
\end{equation}
This situation 
does not seem to have been considered before.
Furthermore, our gluing technique shows that solutions to
\eqref{eq:infinite}
satisfy an \textit{approximate superposition rule}, see
\eqref{eq:superposition} below.
For a finite number of vortex points on $\R^2$, such a rule exists 
formally in the physics literature, and has been rigorously derived in
\cite{Stuart}.
In view of the representation \eqref{eq:superposition},
we can easily analyze the asymptotic behavior of solutions to
\eqref{eq:infinite}
as $\delta\to0^+$, thus obtaining more direct proofs for the
asymptotics derived
in \cite{HongJostStruwe,WangYang}, in the special case
\eqref{eq:periodic}.
\par
In order to state our results, we denote by $\Un$ the unique 
radial solution for the problem:
\begin{equation}
\label{eq:02}
\begin{cases}
-\Delta\Un =1-\text{e}^{\Un}-4\pi N\delta_{0}&\text{in}\ \R^{2}\\
\Un(x) \to 0 &\text{as}\ |x|\to+\infty.
\end{cases}
\end{equation}
Existence, uniqueness and exponential decay as $|x|\to+\infty$
for $\Un$ are established in \cite{JaffeTaubes}, see Section~\ref{sec:singlevortex}
below.
\par
Our main result is the following
\begin{thm}
\label{thm:infinite}
Let $p_j\in\R^2$, $m_j\in\N$, $j\in\N$ satisfy \eqref{eq:constraints}.
There exists a constant $\delta_1>0$ (dependending on $d$ and
$m$ only) such that
for every $\delta\in(0,\delta_1)$ there exists a solution $\ud$ for
\eqref{eq:infinite}.
If the $p_j$'s are doubly periodically arranged in $\R^2$, then $\ud$
is doubly periodic.
Furthermore, $\ud$ satisfies the approximate superposition rule:
\begin{equation}
\label{eq:superposition}
\ud(x)=\sum_{j\in\N}U_{m_j}\left(\frac{|x-p_j|}{\delta}\right)+\omega_\delta,
\end{equation}
where  the error term $\omega_\delta$ satisfies 
$\|\omega_\delta\|_\infty\le C\e^{-c/\delta}$, for some $c>0$ independent of
$\delta$.
In particular, $u$ satisfies the following properties:
\begin{enumerate}
\item[(i)]
$0\le\e^{\ud}<1$, $e^{\ud}$ vanishes exactly at $p_j$, $j\in\N$; 
\item[(ii)]
For every compact subset $K$ of $\R^2\setminus\cup_{j\in\N}\{p_j\}$
there exist $C,c>0$ such that
$\sup_K(1-\e^{\ud})\le C\e^{-c/\delta}$   
as $\delta\to 0^{+}$;
\item[(iii)]
$\delta^{-2}(1-\e^{\ud})\to4\pi\sum_{j\in\N}m_j\delta_{p_j}$
in the sense of distributions, as $\delta\to 0^{+}$.
\end{enumerate}
\end{thm}
We note that $\delta^{-2}(1-\e^{\ud})=2|F_{12}|$.
\par
An outline of this note is as follows.
Our starting point in proving Theorem~\ref{thm:infinite} is to
consider $\delta$ as a scaling parameter.
Setting $\hat u(x)=u(\delta x)$, we have that $\hat u$ satisfies:
\begin{equation}
\label{eq:01}
-\Delta\hat u=1-\e^{\hat u}-4\pi\sum_{j\in\N}m_j\delta_{\hat p_j}
\qquad\mathrm{in\ }\R^2,
\end{equation}
where $\hat p_j=p_j/\delta$. Note that the vortex points $\hat p_j$
``separate" as $\delta\to0^+$.
Section~\ref{sec:singlevortex} contains the necessary properties
of the radial solutions $\Un$ to \eqref{eq:02}. We rely on the results of 
Taubes~\cite{Taubes} for the existence and uniqueness of $\Un$, as well as for 
the exponential decay properties at infinity. We also prove a
necessary non-degeneracy property of $\Un$.
The exponential decay of solutions justifies the following approximate superposition
picture 
for small values of $\delta$, i.e., for vortex points $\hat p_j$ which are ``far
apart":
\begin{equation}
\label{eq:introsuperposition}
\hat u(x)\approx\sum_{j\in\N}U_{m_j}\left(|x-\hat p_j|\right).
\end{equation}
In fact, we take the following preliminary form
of the superposition rule:
\begin{equation}
\label{eq:introansatz}
\hat{u}=\sum_{j\in\N}\phij U_{m_j}(x-\hat p_j)+z,
\end{equation}
as an \textit{ansatz} for $\hud$. 
Here, radial solutions centered at $\hat p_j$
are ``glued" together by the functions $\phij$, which belong to
a suitable locally finite partition of unity. 
Section~\ref{sec:partition}
contains the definition and the main properties of the partition,
as well as of the appropriate functional spaces $\Xd,\Yd$,
which are also obtained by ``gluing" $H^1(\R^2)$ and $L^2(\R^2)$,
respectively.
Hence, we are reduced to show that for small values of $\delta$
there exists an exponentially small ``error" $z$ such that
$\hat{u}$ defined by \eqref{eq:introansatz} is a solution for \eqref{eq:01}.
The existence of such a $z\in\Xd$ is the aim of Section~\ref{sec:shadowing}
(see Proposition~\ref{prop:z}). To this end we use the shadowing lemma.
We characterize $z$ by the property
$\Fd(z)=0$, where $\Fd:\Xd\to\Yd$ is suitably defined. 
The non-degeneracy property of $\Un$ is essential in order to
prove that the operator $D\Fd(0)$ is invertible, and that its inverse is bounded
independently of $\delta>0$ (Lemma~\ref{lem:op-inv}).
At this point, the Banach fixed point argument
applied to $\mathbb I-(D\Fd(0))^{-1}\,\Fd$ yields the existence of
the desired error term $z$.
In Section~\ref{sec:proof} we show that periodically arranged vortex points
lead to periodic solutions, that \eqref{eq:introansatz} implies
\eqref{eq:superposition} and we derive the asymptotic behavior of solutions,
thus concluding the proof of Theorem~\ref{thm:infinite}.
For the reader's convenience, following the monograph 
of Jaffe and Taubes~\cite{JaffeTaubes}, we outline in an appendix
the derivation of equation~\eqref{eq:Taubes} for smooth, finite action
critical points to the action of $\Ed$ on $\R^2$, 
as well as some properties of solutions
to \eqref{eq:Taubes}, which imply the necessary properties of $\Un$.
\par
Although we have chosen
to consider the abelian Higgs model for the sake of simplicity,
it will be clear from the proof that our method may be 
adapted to many other self-dual gauge theories
as considered, e.g., in the monograph \cite{Yang}.
\par
Henceforth, unless otherwise stated, we denote by $C,c>0$ general
constants 
independent of $\delta>0$ and of $j\in\N$.
\section{Single vortex point solutions}
\label{sec:singlevortex}
In this section we consider the solution
$\Un$ to the radially symmetric equation \eqref{eq:02}.
We refer to \cite{JaffeTaubes,Taubes} for the proof 
of the existence and uniqueness of  $\Un$ (see also the Appendix).
We collect in the following lemma some properties of $\Un$ that will
be needed in the sequel. For every $r>0$, we denote 
$B_{r}=\{x\in\R^{2}\,:\,|x|<r\}$.
\begin{lemma}
\label{lem:Ur}
The following properties hold:
\begin{enumerate}
\item[(i)] 
$ \e^{\Un (x)}<1$
for any $x\in\R^2$.
\item[(ii)] 
For every $r>0$ there exist constants $C_N>0$ and  $\alpha_N>0$ 
depending on $r$ and $N$ such that
\begin{equation*}
\label{eq:expdecay}
|1-\e^{\Un(x)}|+|\nabla\Un(x)|+|\Un(x)|\le C_N\e^{-\alpha_N|x|},
\end{equation*}
for all $x\in\R^2\setminus B_r$.
\end{enumerate}
\end{lemma} 
\begin{proof}
Property (i) follows  by the maximum principle.
In order to establish (ii), we note that the 
estimate $|1-\e^{\Un(x)}|\le C_N\e^{-\beta_N|x|}$ for some $\beta_N>0$
depending on $N$
was established by Taubes (\cite{Taubes}, Theorem~III.1.1),
see the Appendix.
In view of (i), it follows that for all $|x|\ge r$
we have
\[
|\Un(x)|=\frac{|\Un(x)|}{1-\e^{\Un(x)}}\left(1-\e^{\Un(x)}\right)
\le C\e^{-\beta|x|}. 
\]
In order to estimate the decay of $|\nabla\Un|$,
we set $A=B_{4r}\setminus \overline B_r$, and for all $R\ge r$ we define 
$\Ar=B_{4Rr}\setminus \overline B_{Rr}$, $\Ar'=B_{3Rr}\setminus 
\overline B_{2Rr}$.
For $y\in A$, we consider $u_{R}(y)=\Un(Ry)$. Then $u_{R}$ satisfies
$-\Delta u_{R}=f_{R}$ in $A$ with $f_R$ given by
$f_{R}(y)=R^2(1-\exp\{\Un(Ry)\})$.
We recall the standard elliptic estimate for $u_{R}$ 
(see, e.g., \cite{Gilbarg-Trudinger} Theorem~3.9):
\begin{equation*}
\sup_A\disty|\nabla u_{R}(y)|
\le C\left(\sup_A|u_{R}|+\sup_A\disty^2|f_{R}(y)|\right),  
\end{equation*}
where $\disty=\mathrm{dist}(y,\partial A)$ and $C>0$ is
independent of $R$.
In terms of $\Un$, the above estimate yields
\begin{equation}
\sup_{\Ar}\distx|\nabla\Un(x)|
\le
C\left(\sup_{\Ar}|\Un|+\sup_{\Ar}\distx^2(1-\exp\{\Un(x)\})\right).
\end{equation}
where $\distx=\mathrm{dist}(x,\partial\Ar)=R\disty$.
Hence, we have for any $x \in \Ar'$
\begin{equation}
|\nabla\Un(x)| \leq \sup_{\Ar'}\frac{\distx}{R}|\nabla\Un(x)| \leq C
R\e^{-\beta R} \leq C|x|\e^{-\frac{\beta}{3}|x|}
\end{equation}
and we conclude that
\begin{equation}
|\nabla\Un(x)| \leq C\e^{-\alpha|x|}\qquad\qquad\forall|x|\geq r
\end{equation}
for some constant $\alpha>0$.
\end{proof}
We consider the bounded linear operator 
\begin{equation*}
\Ln=-\Delta+\e^{\Un}:H^2(\R^{2})\to L^2(\R^{2}).
\end{equation*}
It is known \cite{JaffeTaubes, Taubes} that $\Un$ corresponds to 
the unique minimum of a strictly convex functional,
and therefore it is the unique solution to \eqref{eq:02}, see the Appendix.
In order to apply the shadowing lemma,
we further have to show that $\Un$ is non-degenerate,
in the sense of the following
\begin{lemma}
\label{lem:inv}
The operator $\Ln$ is invertible and for every $N>0$ there exists
$C_N>0$
such that $\|\Ln^{-1}\|\le C_N$. 
\end{lemma}
\begin{proof}
It is readily seen that $\Ln$ is injective. Indeed, suppose
$\Ln u=0$ for some $u\in H^2(\R^{2})$. Multiplying by $u$ and
integrating on $\R^2$
we have:
\[
\int|\nabla u|^2+\int\e^{\Un}u^2=0.
\]
Therefore, $u=0$.
Now we claim that $\Ln$ is a Fredholm operator. Indeed,
we write
\[
\Ln=(-\Delta+1)(\mathbb{I}-T),
\]
with $T=(-\Delta+1)^{-1}(1-\e^{\Un}):H^2(\R^{2})\to H^2(\R^{2})$. 
Clearly, $T$ is continuous.
Let us check that $T$ is compact.
To this end, let $\un\in H^2(\R^{2})$, $\|\un\|_{H^2}=1$. We have to
show that
$T\,\un$ has a convergent subsequence. Note that by the Sobolev
embedding 
\begin{equation}
\label{eq:Sobolev}
\|u\|_{L^\infty(\R^2)}\le C_S\|u\|_{H^2(\R^2)},
\end{equation}
for all $u\in H^2(\R^2)$, we have
$\|\un\|_\infty\le C'$,
for some $C'>0$ independent of $n$, and there exists $\ui$,
$\|\ui\|_{H^2}\le1$, such that
$\unk\to\ui$ strongly in $L_{\mathrm{loc}}^2$ for a subsequence
$\unk$.
Now, by Lemma~\ref{lem:Ur}, for any fixed $\eps>0$, there exists
$R>0$ such that
$\|1-\e^{\Un}\|_{L^2(\R^2\setminus B_R)}\le\eps$. Consequently,
$\|(1-\e^{\Un})(\unk-\ui)\|_{L^2(\R^2\setminus B_R)}\le 2C'\eps$.
On the other hand, $\|(1-\e^{\Un})(\unk-\ui)\|_{L^2(B_R)}\to0$.
We conclude that $(1-\e^{\Un})(\unk-\ui)\to0$ in $L^2$.
In turn, we have
$T(\unk-\ui)=(-\Delta+1)^{-1}(1-\e^{\Un})(\unk-\ui)\to0$
in $H^2$, which implies that $T$ is compact.
It follows that $\Ln$ is a Fredholm operator. Consequently, $\Ln$ is
also surjective. 
At this point, the Open Mapping Theorem concludes the proof.
\end{proof}
\section{A partition of unity}
\label{sec:partition}
In this section we introduce a partition of unity and we prove some
technical results
which will be needed in the sequel.
Let $p_j\in\R^2$, $j\in\N$ be the vortex points.
By assumption \eqref{eq:constraints},  $\ro=d/8= \inf_{j\neq k}|p_j-p_k|/8>0$.
We consider the set $K = (-\frac{3}{4} \ro, \frac{3}{4} \ro )  \times
(-\frac{3}{4} \ro, \frac{3}{4} \ro )$. Then for any $\underline n
\in \Z^2$, we introduce $K_{\underline n}=K+\underline n \ro $. The
collection of sets 
$\{ K_{\underline n} \}_{\underline n\in \Z^{2}}$ is a locally finite covering of $\R^{2}$. 
We consider an associated partition of unity defined as follows: let
$0 \leq \phi \in C^{\infty}_{c}(K )$ be such that $\sum_{\underline n
\in \Z^{2}} \phi_{\underline n} (x) = 1$ pointwise, where
$\phi_{\underline n} (x) = \phi(x -  \underline n \ro)$.
Then, for any $j \in \N$, we introduce the set 
$$
N_{j}=\{\underline n\in\Z^2\,:\,d(p_{j},K_{\underline n})
< \frac{1}{4}\ro\}, 
$$
note that  the cardinality of $N_{j}$ is uniformly bounded, namely
$|N_{j}| \leq 4 $ for any $j \in \N$.
For any $j \in \N$, we set
$$
B_{j}=\bigcup_{ \underline n \in N_{j}} K_{ \underline n},  \qquad
\varphi_j (x)= \sum_{ \underline n \in N_{j}} \phi_{\underline n}(x). 
$$
Let $\mathcal{I}:\N\to\Z^2\setminus\bigcup_{j \in \N} N_j$
be a bijection. We set
$$
Q_j = K_{ \mathcal{I} (j)}, \qquad \psi_j(x)= \phi_{\mathcal{I}
(j)}(x).
$$
Then $\{B_j,Q_j\}_{j\in\N}$ is a locally finite open covering of
$\R^2$
with the property that $B_j\cap B_k=\emptyset$ for every $k\neq j$.
Moreover $\{\varphi_j,\psi_j\}$ is a   partition of unity associated
to
$\{B_j,Q_j\}_{j\in\N}$, such that
\[\mathrm{supp} \, \varphi_j\subset B_j,\qquad\mathrm{supp}
\,\psi_j\subset Q_j,\]
and such that
\begin{align*}
\sup_{j\in\N}\{\|\nabla\varphi_j\|_\infty,\ 
\|\nabla\psi_j\|_\infty\}<+\infty,
\qquad
\sup_{j\in\N}\{\|D^2\varphi_j\|_\infty,\ 
\|D^2\psi_j\|_\infty\}<+\infty.
\end{align*}
In particular,
\[
0\le\varphi_j,\psi_j\le1
\qquad\mathrm{and}\qquad
\sum_{j\in\N}(\varphi_j(x)+\psi_j(x))= \sum_{\underline n \in \Z^2}
\phi_{\underline n}(x) =1.
\]
For every $j\in\N$, we define a rescaled covering:
\[
\Bj=B_j/\delta,
\qquad\Qj=Q_j/\delta.
\]
Then $\{\phij,\psij\}_{j\in\N}$ defined by
\[
\phij(x)=\varphi_j(\delta x),\qquad\psij(x)=\psi_j(\delta x)
\]
is a partition of unity associated to $\{\Bj,\Qj\}$.
It will also be convenient to define the sets
\begin{equation*}
\Cj=\{x\in\Bj\,:\,\phij(x)=1\}
\qquad j\in\N.
\end{equation*}
Note that 
\begin{align*}
\mathrm{supp}\{\nabla\phij,D^2\phij\}\subset\Bj\setminus\Cj
\end{align*}
and
\begin{align}
\label{eq:gradphidecay}
\sup_{\R^{2}}\{|\nabla\phij|+|\nabla\psij|\}\leq C\delta,
\qquad\sup_{\R^{2}}\{|D^2\phij|+|D^2\psij|\}\leq C\delta^{2}.
\end{align}
For every fixed $x\in\R^2$ we define the following subsets of $\N$:
\begin{equation}
\label{eq:JKdef}
J(x)=\{j\in\N\,:\,\phij(x)\neq0\},
\qquad
K(x)=\{k\in\N\,:\,\hat\psi_k(x)\neq0\}.
\end{equation}
Note that
\begin{equation}
\label{eq:supJK}
\sup_{x\in\R^2}\{|J(x)|+|K(x)|\}<+\infty,
\end{equation}
where $|J(x)|$, $|K(x)|$ denote the cardinality of
$J(x)$, $K(x)$, respectively.
We shall use the following Banach spaces:
\begin{equation*}
\begin{split}
\Xd=&\{u\in H^{2}_{\mathrm{loc}}(\R^{2})\,:\,\, 
\sup_{j\in\N}\{\|\phij u\|_{H^{2}(\R^{2})},\|\psij
u\|_{H^{2}(\R^{2})}\}<+\infty\}\\
\Yd=&\{f\in L^{2}_{\mathrm{loc}}(\R^{2})\,:\,\,
\sup_{j\in\N}\{\|\phij f\|_{L^{2}(\R^{2})},\|\psij
f\|_{L^{2}(\R^{2})}\}<+\infty\}.
\end{split}
\end{equation*}
We collect in the following lemma some estimates that
will be used in the sequel.
\begin{lemma}
\label{lem:basicestimates} 
There exists a constant  $C>0$ such that
for any $u\in\Xd$ and $j\in\N$ we have 
\begin{enumerate}
\item[(i)]
$\|u\|_{H^2(\Bj)}\le C\|u\|_{\Xd}$
\item[(ii)]
$\|u\|_{L^\infty(\R^2)}\le C\|u\|_{\Xd}$.
\end{enumerate}
\end{lemma}
\begin{proof}
(i)
For every fixed $k\in\N$, let
$\mathcal J(k)=\{j\in\N\,:\,\mathrm{supp} \,\psij\cap\mathrm{supp}
\,\hat\varphi_k\neq\emptyset\}$.
Then $\sup_{k\in\N}|\mathcal J(k)|< +\infty$, and 
we estimate:
\begin{align*}
\|u\|_{H^2(\Bj)}=&\|\phij u+\sum_{k\in \mathcal J(j)}\psik
u\|_{H^2(\Bj)}
\le\|\phij u\|_{H^2(\Bj)}+\sum_{k\in \mathcal J(j)}\|\psik
u\|_{H^2(\Bj)}\\
\le&\left(1+|\mathcal J(j)|\right)\|u\|_{\Xd}\le C\|u\|_{\Xd}.
\end{align*}
\par
(ii)
For any fixed $x\in\R^2$ we have, in view of \eqref{eq:Sobolev}
and \eqref{eq:supJK}:
\begin{align*}
|u(x)|=&\sum_{j\in\N}\phij(x)|u(x)|+\sum_{j\in\N}\psij(x)|u(x)|\\
=&\sum_{j\in J(x)}\phij(x)|u(x)|+\sum_{j\in K(x)}\psij(x)|u(x)|\\
\le&\sum_{j\in J(x)}C_S\|\phij u\|_{H^2(\R^2)}
+\sum_{j\in K(x)}C_S\|\psij u\|_{H^2(\R^2)}\\
\le&\sup_{x\in\R^2}(|J(x)|+|K(x)|)C_S\|u\|_{\Xd}=C\|u\|_{\Xd}.
\end{align*}
Hence, (ii) is established.
\end{proof}
We shall also need the following family of functions: 
\[
\gj=\frac{\phij}{\left(\sum_{k\in\N}
(\hat{\varphi}_{k}^{2}+\hat\psi_k^2)\right)^{1/2}},
\qquad
\hj=\frac{\psij}{\left(\sum_{k\in\N}
(\hat{\varphi}_{k}^{2}+\hat\psi_k^2)\right)^{1/2}}.
\]
In view of \eqref{eq:gradphidecay}, it is readily checked that
\begin{lemma}
 \label{lem:phig}
The family $\{\gj, \hj \}_{j\in \N}$ satisfies
$\mathrm{supp}\gj\subset\Bj$, $\mathrm{supp}\hj\subset\Qj$
    and furthermore:
\begin{align}
  \label{eq:sumgh}
   &\sum_{j\in\N}(\gj^{2}+\hj^2)\equiv1\\
    &C^{-1}\phij\le\gj\le C\phij,
&&C^{-1}\psij\le\hj\le C\psij\\
\label{eq:gradphidecay1}
    &\sup_{\R^{2}}\{|\nabla\gj|+|\nabla\hj|\}\leq C\delta,
    &&\sup_{\R^{2}}\{|D^2\gj|+|D^2\hj|\}\leq C\delta^{2}.
    \end{align}
 \end{lemma}  

\section{The shadowing lemma}
\label{sec:shadowing}
For every $j\in\N$ we define
\[
\Uj(x)=U_{m_j}(x-\hat p_j).
\]
We make the following ansatz for solutions $\hat u$ to 
equation \eqref{eq:01}:
\begin{equation}
\label{eq:ansatz}
\hat{u}=\sum_{j\in\N}\phij\hat{U}_{j}+z.
\end{equation}
Our aim in this section is to prove:
\begin{prop}
\label{prop:z}
There exists $\delta_1>0$ such that for all $\delta\in(0,\delta_1)$
there exists $\zd\in\Xd$, such that $\hat u_\delta$ defined by
$\hat u_\delta=\sum_j\phij\hat{U}_{j}+\zd$ is a solution to
\eqref{eq:01}.
Moreover, $\|\zd\|_{\Xd}\le C\e^{-c/\delta}$.
\end{prop}
We note that the functional $\Fd:\Xd\to\Yd$ given by
\[
\Fd(z)=-\Delta
z+\sum_{j\in\N}\phij(1-\e^{\Uj})-(1-\e^{\sum_{j\in\N}{\phij\Uj+z}})
-\sum_{j\in\N}[\phij, \Delta]\Uj
\]
is well-defined and $C^1$. 
Here $[\Delta,\phij]=\Delta\phij+2\nabla\phij\nabla$.
Moreover, if $z\in\Xd$ satisfies $\Fd(z)=0$,
then $\hat u$ defined by \eqref{eq:ansatz} is a solution to
\eqref{eq:01}.
\begin{lemma}
\label{lem:stimaF} For $\delta >0$
sufficiently small, we  have
\begin{equation}
\|\Fd(0)\|_{\Yd}\leq C\e^{-c/\delta}
\qquad\mathrm{as}\,\,\delta \to 0^{+}
\end{equation}
for some constants $C,c>0$ independent of $\delta$.
\end{lemma}
\begin{proof}
Let  
\begin{equation*}
\begin{split}
\mathcal R &
=\sum_{j\in\N}\phij(1-\e^{\Uj})-(1-\e^{\sum_{j\in\N}{\phij\Uj}})
\\
\mathcal C & =\sum_{j\in\N} [\phij, \Delta]\Uj \\
\end{split}
\end{equation*}
Note that $ \{ \mathrm{supp} \, \mathcal R,  \, \mathrm{supp} \,
\mathcal C
\} \subset\cup_{j\in\N}\Bj\setminus\Cj$.
We fix $x\in\cup_j\Bj$. We estimate:
\begin{align*}
|\mathcal R(x)|
\le& \sup_{j \in \N}\|\phij(1-\e^{\Uj})\|_{L^\infty(\Bj\setminus\Cj)}
+ \sup_{j \in \N} \|1-\e^{\phij\Uj}\|_{L^\infty(\Bj\setminus\Cj)}\\
\le&C \,\sup_{j\in \N}\|\Uj\|_{L^\infty(\Bj\setminus\Cj)}
\le C_1\e^{-c_1/\delta}.
\end{align*}
On the other hand, in view of \eqref{eq:gradphidecay} and Lemma~\ref{lem:Ur}, 
for $x \in
\cup_j\Bj$, we have
\begin{equation}
   \begin{split}
\label{eq:Ci-1 }
|\mathcal C(x)|
\le&\sup_{j \in \N} \|\,[\Delta,\phij]\Uj
\|_{L^\infty(\Bj\setminus\Cj)}
\\
\le&C\,(\sup_{j\in \N}
\|\Uj\,\Delta\phij\|_{L^\infty(\Bj\setminus\Cj)} +
\sup_{j\in \N}
\|\,|\nabla\Uj|\,|\nabla\phij|\,\|_{L^\infty(\Bj\setminus\Cj)})
\le
C_{2}\text{e}^{-c_{2}/\delta} \\
\end{split} 
\end{equation}
for some positive constants $c_{2}, C_{2}>0$ independent of $\delta$.
Hence, we conclude that, as $\delta \to 0^{+}$,
\begin{equation}
   \|\Fd(0)\|_{\Yd}\leq C \sup_{j \in \N}  (\|\mathcal R
   \|_{L^{2}(\Bj)}
   + \|\mathcal C \|_{L^{2}(\Bj)} ) \le  C\text{e}^{-c/\delta}
   \end{equation}
   for some constants $C,c>0$ independent of $\delta>0$.
\end{proof}
Now, we consider the operator $\Ld\equiv D\Fd(0 ):\Xd\to\Yd$
given by
\begin{equation*}
\Ld=-\Delta+\e^{\sum_{j\in\N}\phij\Uj}.
\end{equation*}
For every $j\in\N$, we define the operators:
\[
\Lj=-\Delta+\e^{\Uj}.
\]
It will also be convenient to define:
\[
\Lo=-\Delta+1.
\]
We readily check the that the following holds:
\begin{lemma}
\label{lem:LdLj}
There exists a constant $C>0$ such that for any $u\in \Xd$ and $j \in
\N$ we have
\begin{align*}
\|(\Ld-\Lj)\phij u\|_{L^2}\le C\e^{-c/\delta}\|\phij h\|_{L^2},\\
\|(\Ld-\Lo)\psij u\|_{L^2}\le C\e^{-c/\delta}\|\psij h\|_{L^2}. \\
\end{align*}
\end{lemma}
\begin{proof}
For any $j \in \N$, by Lemma \ref{lem:Ur}, we have as 
$\delta \to 0^{+}$, 
\begin{equation}
\begin{split}
\label{eq:stima2}
\| (L_{\delta}-\Lj)\phij u\|_{L^{2}}&
\leq (\|1- \e^{\Uj}\|_{L^{\infty}(\Bj\setminus\Cj)} 
+\|1-\e^{\phij \Uj}\|_{L^{\infty}(\Bj\setminus\Cj)} )
\|\phij u \|_{L^2}\\
\leq & C\|1- \e^{\Uj}\|_{L^{\infty}(\Bj\setminus\Cj)} \|\phij u\|_{L^2}
\leq C\e^{-c/\delta}\|\phij u\|_{L^2}.
\end{split}
\end{equation}
Let $\mathcal K(j)=\{k\in\N\,:\,\mathrm{supp}\phik\cap\mathrm{supp}\psij
\neq\emptyset\}$. Then $\sup_{j\in\N}|\mathcal K(j)|<+\infty$
and we estimate, as 
$\delta \to 0^{+}$,
   \begin{equation}
   \begin{split}
   \label{eq:stima4}
\|(L_{\delta}-\hat{L}_{0})\hat{\psi}_{j}u\|_{L^{2}} \leq &
\|(1-\e^{\sum_{k\in\mathcal K(j)}\phik \Uk} )  \hat{\psi}_{j}u\|_{L^{2}}\\
  \leq &\sup_{k \in \mathcal K(j)} \| 1- \e^{\phik\Uk}  
  \|_{L^{\infty}(\Bk\setminus\Ck)} \|\hat{\psi}_{j}u\|_{L^{2}} 
 \leq C\e^{-c/\delta}\|\psij u\|_{L^2}. 
   \end{split}
   \end{equation}
\end{proof}
Now we prove an essential non-degeneracy property of $\Ld$:
\begin{lemma}
\label{lem:op-inv} 
There exists $\delta_{0} >0$ such that for any 
$\delta\in (0,\delta_{0})$,  the operator $\Ld$ is invertible. 
Moreover, $\Ld^{-1}:\Yd\to\Xd$ is uniformly bounded with respect to
$\delta\in(0, \delta_{0})$.
\end{lemma}
\begin{proof}
Following a gluing technique introduced in \cite{Angenent},
we construct an ``approximate inverse" $\Sd:\Yd\to\Xd$ for $\Ld^{-1}$
as follows:  
\begin{equation}
\Sd=\sum_{j\in\N}\left(\gj\Lj^{-1}\gj+\hj\Lo^{-1}\hj\right),
\end{equation}
where $\gj,\hj$ are the functions introduced in Section~\ref{sec:partition}.
We claim that the operator $\Sd$ is well-defined and uniformly bounded
with respect to $\delta$. That is, we claim that
\begin{equation}
\label{eq:Sd}
\|\Sd f\|_{\Xd}\le C\|f\|_{\Yd}
\end{equation} 
for some $C>0$ independent of $f\in\Xd$ and of $\delta>0$.
\par
Indeed, for any $f\in\Yd$ we have
\begin{align*}
\|\Sd f\|_{\Xd}
=\sup_{k\in\N}\{\,\|\hat\varphi_k\sum_{j\in\N}(\gj\Lj^{-1}\gj 
&+\hj\Lo^{-1}\hj)f\|_{H^2},\\
&\|\hat\psi_k\sum_{j\in\N}(\gj\Lj^{-1}\gj
+\hj\Lo^{-1}\hj)f\|_{H^2}\}.
\end{align*}
We estimate, recalling the properties of $\phij$ and $\gj$:
\begin{align*}
\|\phik\sum_{j\in\N}&\gj\Lj^{-1}\gj f\|_{H^2}
=\|\phik\gk\Lk^{-1}\gk f\|_{H^2}\\
\le&C\|\Lk^{-1}\gk f\|_{H^2}\le C\|\gk f\|_{L^2}
\le C\|\phik f\|_{L^2}\le C\|f\|_{\Yd}.
\end{align*}
We have:
\begin{align*}
\|\phik\sum_{j\in\N}\hj\Lo^{-1}\hj f\|_{H^2}
\le\|\phik\sum_{j\in\mathcal J(k)}\hj\Lo^{-1}\hj f\|_{H^2}
\le\sum_{j\in\mathcal J(k)}\|\phik\hj\Lo^{-1}\hj f\|_{H^2},
\end{align*}
where
$\mathcal J(k)=\{j\in\N\,:\,\mathrm{supp} \,\psij\cap\mathrm{supp}
\,\hat\varphi_k\neq\emptyset\}$
satisfies $\sup_{k\in\N}|\mathcal J(k)|< +\infty$. 
In view of Lemma \ref{lem:phig} and 
Lemma~\ref{lem:inv}, we estimate:
\begin{align*}
\sum_{j\in\mathcal J(k)}&\|\phik\hj\Lo^{-1}\hj f\|_{H^2}
\le C\sum_{j\in\mathcal J(k)}\|\Lo^{-1}\hj f\|_{H^2}
\le C\sum_{j\in\mathcal J(k)}\|\hj f\|_{L^2}\\
\le&\sum_{j\in\mathcal J(k)}\|\psij f\|_{L^2}
\le|\mathcal J(k)|\sup_{j\in\N}\|\psij f\|_{L^2}
\le C\|f\|_{\Yd}.
\end{align*}
Therefore,
\[
\sup_{k\in\N}\|\phik\sum_{j\in\N}\hj\Lo^{-1}\hj f\|_{H^2}
\le C\|f\|_{\Yd}.
\]
Similarly, we obtain that
\begin{align*}
&\sup_{k\in\N}\|\psik\sum_{j\in\N}\gj\Lj^{-1}\gj f\|_{H^2}\le
C\|f\|_{\Yd},
&&\sup_{k\in\N}\|\psik\sum_{j\in\N}\hj\Lo^{-1}\hj f\|_{H^2}\le
C\|f\|_{\Yd}.
\end{align*}
and \eqref{eq:Sd} follows.
\par
Now, we claim that there
exists  $\delta_0$ such that for any $\delta\in (0,\delta_0)$, 
the operator $\Sd\Ld:\Xd\to\Xd$
is invertible, and furthermore $\|\Sd\Ld\|\le C$ for some $C>0$
independent of $\delta>0$. 
We note that $(\Ld-\Lj)\gj : \Xd \to \Yd$ and  
$(\Ld-\Lo)\hj : \Xd \to \Yd$ are  well-defined bounded linear
operators.  
Thus, recalling \eqref{eq:sumgh} we decompose:
\begin{equation}
\begin{split}
\label{eq:SA}
\Sd\Ld=&\I_{\Xd}+\sum_{j\in\N}\gj\Lj^{-1}(\gj L_{\delta}-\Lj\gj)
+\sum_{j\in\N}\hj\Lo^{-1}(\hj L_{\delta}-\Lo\hj)
\\
=&\I_{\Xd}+\sum_{j\in\N}\gj\Lj^{-1}(L_{\delta}-\Lj)\gj
+\sum_{j\in\N}\hj\Lo^{-1}(L_{\delta}-\Lo)\hj
+\sum_{j\in\N}\gj\Lj^{-1}[\Delta,\gj]\\
&\qquad\qquad+\sum_{j\in\N}\hj\Lo^{-1}[\Delta,\hj].
\end{split}
\end{equation}
Hence, it suffices to prove that the last four terms in
\eqref{eq:SA} are sufficiently small, in the operator norm, provided
$\delta>0$ is sufficiently small. 
By Lemma~\ref{lem:LdLj} and Lemma~\ref{lem:phig} we have, for any $u\in\Xd$, 
\begin{equation*}
    \begin{split}
    \|\sum_{j\in\N}\gj\Lj^{-1}&( L_{\delta}-\Lj)\gj u \|_{\Xd} \\
     = \sup_{k\in\N} &\{\|\phik\sum_{j\in\N}\gj\Lj^{-1}(
L_{\delta}-\Lj)\gj
    u\|_{H^{2}}, 
    \|\psik\sum_{j\in\N}\gj\Lj^{-1}( L_{\delta}-\Lj)\gj u\|_{H^{2}}\}\\
    \leq &
    C\sup_{k \in \N} \|\Lj^{-1}( L_{\delta}-\Lk)\gk u \|_{H^{2}} \leq 
    C\sup_{k \in \N} \|( L_{\delta}-\Lk)\gk u \|_{L^{2}}\\
\leq &C
    \text{e}^{-c/\delta} \sup_{k \in \N}\|\phik u\|_{L^{2}}
    \leq  C
    \text{e}^{-c/\delta}\| u\|_{\Xd}.
    \end{split}
    \end{equation*}
Similarly, for $u\in\Xd$, we have:
\begin{equation*}
    \begin{split}
\|\sum_{j\in\N}\gj\Lj^{-1} &[\Delta,\gj] u \|_{\Xd}\\
= &\sup_{k\in\N}\{\|\phik\sum_{j\in\N}\gj\Lj^{-1}[\Delta,\gj]
u\|_{H^{2}},
\|\psik\sum_{j\in\N}\gj\Lj^{-1}[\Delta,\gj] u\|_{H^{2}}\}\\
\leq &C\sup_{k\in\N}\|\Lk^{-1}[\Delta,\gk] u\|_{H^{2}}
\le C\sup_{k\in\N}\|[\Delta,\gk] u\|_{L^{2}}.
\end{split}
  \end{equation*}
Recalling that $[\Delta,\gk] u=2\nabla  u\nabla\gk+ u\Delta\gk$,
 by \eqref{eq:gradphidecay1} and Lemma~\ref{lem:basicestimates}--(i) we derive that
\begin{align*}
\|[\Delta,\gk] u\|_{L^{2}}
\le C\delta\| u \|_{H^1(\hat B_k)}\le C\delta\| u\|_{\Xd}.
\end{align*}
The remaining terms are estimated similarly.
Hence, $\|\Sd\Ld-\I_{\Xd}\|\to0$ as $\delta\to0^+$.
Now we observe that
$\Ld^{-1}=(S_{\delta}L_{\delta})^{-1}\Sd$.
It follows that for any $f\in\Yd$ we have
\begin{equation}
\|\Ld^{-1}f\|_{\Xd}=\|(S_{\delta}L_{\delta})^{-1}\Sd f\|_{\Xd} 
\leq C \|\Sd f\|_{\Yd} 
\leq C\|f\|_{\Yd}
\end{equation}
with $C >0$ independent of $\delta$.
Hence, $\Ld$ is invertible and its inverse is bounded independently
of $\delta$, as asserted.
\end{proof}
Now we can provide the 
\begin{proof}[Proof of Proposition~\ref{prop:z}]
We use the Banach fixed point argument.
For any $\delta \in (0, \delta_{0})$,
with $\delta_{0} >0$ given by Lemma~\ref{lem:op-inv}, 
we introduce the nonlinear  map 
$G_{\delta}\in C^{1}(\hat{X}_{\delta},\hat{X}_{\delta})$ defined by
\begin{equation}
\label{eq:mappa}
G_{\delta}(z)=z-L_{\delta}^{-1}F_{\delta}(z).
\end{equation}
and the set 
\begin{equation}
    \mathcal{B}_{R} =\{ u \in \Xd \, : \, \|u\|_{\Xd} \leq R \}
    \end{equation}  
Then, fixed points of $G_{\delta}$ correspond to solutions of
the functional equation $F_{\delta}(z)= 0$.
First, note that   $DG_{\delta}(0) = 0$ and that
\[
DF(z)=-\Delta+\e^{\sum_{j\in\N}\phij\Uj+z}.
\]
By Lemma~\ref{lem:op-inv}, for any $z\in\Xd$ and $u\in\Xd$ we
have
\begin{align*}
\|DG_{\delta}&(z)u\|_{\Xd}=\|(DG_{ \delta }(z)-DG_{\delta}(0))
u\|_{\Xd}
=\|L_{\delta }^{-1}(DF_{\delta}(z)-L_{\delta })u\|_{\Xd}\\
\leq&C\|(DF_{\delta}(z)-L_{\delta })u\|_{\Yd}
=C\|\e^{\sum_{j\in\N}\phij\Uj}(\e^z-1)u\|_{\Yd}
\le C\|(\e^z-1)u\|_{\Yd}.
\end{align*}
By the elementary inequality $\e^t-1\le Ct\e^t$, for all $t>0$,
where $C>0$ does not depend on $t$, and in view of
Lemma~\ref{lem:basicestimates},
we have
\begin{equation}
\|\text{e}^{z}-1\|_\infty\le\e^{\|z\|_\infty}-1
\le C\|z\|_\infty\e^{\|z\|_\infty}
\le C\|z\|_{\Xd}\e^{\|z\|_{\Xd}}.
\end{equation}
Hence, 
\begin{equation*}
\begin{split}
\|DG_{\delta}(z)u\|_{\hat{X}_{\delta}} 
\leq C \|(\e^z-1)u\|_{\hat{Y}_{\delta}}
\le C\| z\|_{\Xd}\e^{\|z\|_{\Xd}}\|u\|_{\Yd}
\le C\|z\|_{\Xd}\e^{\|z\|_{\Xd}}\|u\|_{\Xd}.
\end{split}
\end{equation*}
Consequently, there exists $R_0>0$ such that for every $R\in(0,R_0)$
we have
\begin{equation}
\|DG_{\delta}(z)\|\leq\frac{1}{2} ,\qquad \forall z
\in \mathcal{B}_{R}
\end{equation}
and for all $\delta>0$.
Now,
\begin{equation}
\begin{split}
\|G_{\delta}(z)\|_{\hat{X}_{\delta}}&
\leq\|G_{\delta}(z)-G_{\delta }(0)\|_{\hat{X}_{\delta}} 
+\|G_{\delta}(0)\|_{\hat{X}_{\delta}}\\
&\leq\frac{1}{2}\|z\|_{\hat{X}_{\delta}}+
\|L^{-1}_{\delta}F_{\delta}(0)\|_{\hat{X}_{\delta}}.
\end{split}
\end{equation}
By Lemma~\ref{lem:op-inv} and Lemma~\ref{lem:stimaF}, we
have that:
\begin{equation}
\|L^{-1}_{\delta}F_{\delta}(0)\|_{\hat{X}_{\delta}} \leq
C\|F_{\delta }(0)\|_{\hat{Y}_{\delta}}\leq C_0\text{e}^{-c_0/\delta}.
\end{equation}
Choosing $R=R_\delta=2C_0\text{e}^{-c_0/\delta}$, we obtain that
$G_\delta(B_{R_\delta})\subset B_{R_\delta}$.
Hence, 
$G_{\delta}$ is a strict contraction in
$B_{R_\delta}$,
for any $\delta \in  (0, \delta_{1})$. By the Banach fixed-point
theorem, for any $\delta \in (0, \delta_{1})$, there exists a unique 
$z_{\delta}\in{B}_{R_\delta}$, such that
$F_{\delta}(z_{\delta}) = 0$.
\end{proof} 

\section{Proof of Theorem~\ref{thm:infinite}}
\label{sec:proof}
In this section we finally provide the proof of
Theorem~\ref{thm:infinite}.
In view of Proposition~\ref{prop:z}, the function $\hud$ defined by
\begin{equation}
\label{eq:sol}
\hud=
\sum_{j\in\N}\phij\Uj+\zd
\end{equation}
is a solution to equation~\eqref{eq:01}.
Consequently, $\ud$ defined by
\begin{equation}
\label{eq:sol1}
\ud(x)=\hud\left(\frac x\delta\right)
=\sum_{j\in\N}\varphi_j(x)U_{m_{j}}\left(\frac{x-p_j}{\delta}\right)
+\zd\left(\frac{x}{\delta}\right)
\end{equation}
is a solution to \eqref{eq:infinite}.
Now, we want to prove that if the $p_j$'s  are doubly periodically 
arranged in $\R^2$,
then $\ud$ is in fact  a 
doubly periodic solution to \eqref{eq:periodic}. 
Recall from Section~\ref{sec:intro} that the $p_j$'s  are doubly
periodically 
arranged in $\R^2$ 
if there exists $s\in\N$ such that for any $k\in\N$, $k>s$ there exist
$j\in\{1,2,\ldots,s\}$ and $m,n\in\Z$ such that $p_k=p_j+m\underline
e_1+n\underline e_2$,
where $\underline e_1,\underline e_2$ are the 
unit vectors in $\R^2$.
We define $\ek=\underline e_k/\delta$, $k=1,2$.
Equivalently, we show:
\begin{lemma}
\label{lem:period}
Suppose the vortex points $p_j$,
$j\in\N$, are doubly periodically arranged in $\R^2$.
Then
$\hud(x+\ek)=\hud(x)$  
for any $x\in\R^{2}$ and for $k=1,2$.
\end{lemma}
\begin{proof}
We may assume that $\phij(x+\ek)=\phij(x)$, $\psij(x+\ek)=\psij(x)$,
for any $j\in\N$, $x\in\R^{2}$, $k=1,2$.
Then,
\begin{equation*}
\hud(x+\ek) =
\sum_{j\in\N}\phij(x)\Uj(x)+\zd(x+\ek). 
\end{equation*}
Hence, it is sufficient to prove that $\zd(x+\ek)=\zd(x)$,
for every $x\in\R^2$ and for $k=1,2$.
First, we claim that $\zd(\,\cdot\,+\ek)\in\mathcal{B}_{\Rd}$. Indeed,
for every $j\in\N$ there exists exactly one $j'\in\N$ such that
\begin{equation}
\|\phij\zd(\,\cdot\,+\ek)\|_{H^{2}}=\|\hat\varphi_{j'}\zd\|_{H^{2}}.
\end{equation}
Hence, we obtain
\begin{equation}
\|\zd(\,\cdot\,+\ek)\|_{\Xd} = 
\|\zd\|_{\Xd}\leq\Rd.
\end{equation}
Moreover, if $F_{\delta}(\zd)=0$ we also have 
$F_{\delta}(\zd(\,\cdot\,+\ek))=0$.
Therefore, $\zd(\,\cdot\,+ \ek)$ is a
fixed point of $G_{\delta}$ in $\mathcal{B}_{\Rd}$. 
By uniqueness, we conclude  that $\zd(\,\cdot\,+\ek)=\zd$,
$k=1,2$, as asserted.
\end{proof}
\begin{lemma}
\label{lem:superposition}
The solution $\ud$ defined in \eqref{eq:sol1} satisfies
the approximate superposition rule:
\begin{equation}
\label{eq:deltasuperposition}
\ud(x)=\sum_{j\in\N}U_{m_j}\left(\frac{x-p_j}{\delta}\right)
+\omega_\delta(x),
\end{equation}
with $\|\omega_\delta\|_\infty\le C\e^{-c/\delta}$.
\end{lemma}
\begin{proof}
In view of \eqref{eq:sol1} and of the definition 
of $J(x)$ in Section~\ref{sec:partition}, we have
\begin{equation*}
\ud(x)=\sum_{j\in J(x)}U_{m_j}\left(\frac{x-p_j}{\delta}\right)
+\tilde\omega_\delta(x),
\end{equation*}
where
\[
\tilde\omega_\delta(x)=-\sum_{j\in J(x)}(1-\varphi_j(x))
U_{m_j}\left(\frac{x-p_j}{\delta}\right)+\zd\left(\frac{x}{\delta}\right).
\]
We estimate:
\begin{equation*}
\left\|\sum_{j\in J(x)}(1-\varphi_j(x))
U_{m_j}\left(\frac{x-p_j}{\delta}\right)\right\|_\infty
\le \sum_{j\in J(x)}\sup_{\R^2\setminus C_j}
\left|U_{m_j}\left(\frac{x-p_j}{\delta}\right)\right| \le
C\e^{-c/\delta}.
\end{equation*}
On the other hand, we readily have 
\[
\|\zd\left(\frac{\, \cdot \,}{\delta}\right)\|_\infty
=\|\zd\|_\infty\le C\e^{-c/\delta}.
\]
Therefore, $\|\tilde\omega_\delta\|_\infty\le C\e^{-c/\delta}$.
We have to show that
\begin{equation}
\label{eq:remainder}
\left\|\sum_{j\not\in
J(x)}U_{m_j}\left(\frac{x-p_j}{\delta}\right)\right\|_\infty
\le C\e^{-c/\delta}.
\end{equation}
To this end, we fix $x\in\R^2$ and for every $N\in\N$ we define 
$B_N=\{y\in\R^{2}\,:\,|y -x |<r_{0}N\}$.
Then,
\begin{equation*}
\sum_{j\not\in J(x)}U_{m_j}\left(\frac{x-p_j}{\delta}\right)
=\sum_{N\in\N}\sum_{p_j\in\overline{B_{N+1}}\setminus B_N}
U_{m_j}\left(\frac{x-p_j}{\delta}\right)
\end{equation*}
Since $\inf_{j \not=k} |p_{j}-p_{k} | > r_{0}$ there exists $C>0$ 
independent of $N\in\N$ and
of $x\in\R^2$ such that
\begin{equation}
\label{eq:N}
\left|\{p_j\in\overline{B_{N+1}}\setminus B_N\}\right|\le C N.
\end{equation}
Hence, we estimate:
\begin{align*}
\left|\sum_{j\not\in
J(x)}U_{m_j}\left(\frac{x-p_j}{\delta}\right)\right|
\le C\sum_{N\in\N}N\e^{-cN/\delta}\le C\e^{-c/\delta}.
\end{align*}
This implies \eqref{eq:deltasuperposition}.
\end{proof}
We are left to analyze the asymptotic
behavior of $\ud$ as $\delta\to 0^{+}$.  
Such a behavior is a straightforward consequence of \eqref{eq:sol1}.
\begin{lemma}
\label{lem:asymp}
Let 
$u_{\delta}$ be  given by \eqref{eq:sol1}. The
following properties hold:
\begin{enumerate}
\item[(i)]
$\text{e}^{u_{\delta}} < 1$ on $\R^2$ and vanishes exactly at  $p_{j}$
with multiplicity $2m_{j}$, $j\in\N$;
\item[(ii)]
For every compact subset $K$ of $\R^2\setminus\cup_{j\in\N}\{p_j\}$
there exist $C,c>0$ such that
$1-\text{e}^{u_{\delta}}\le C\e^{-c/\delta}$ as
$\delta \to 0^{+}$;
\item[(iii)]
$\delta^{-2}(1-\text{e}^{u_{\delta}})\to4 \pi \sum_{j\in\N}m_{j}
\delta_{p_{j}}$ in the sense of distributions, as $\delta \to
0^{+}$.
\end{enumerate}
\end{lemma}
\begin{proof}
(i) Since $u_{\delta} $ is a solution of equation
\eqref{eq:infinite},
 $ \text{e}^{u_{\delta}}< 1$
follows by the maximum principle.
Moreover, since 
\begin{equation}
\label{eq:L1}
U_{m_{j}}((x-p_{j})/\delta)= 
\ln|x-p_{j}|^{ 2m_{j}}+ v_{j} 
\end{equation}
with $v_{j}$ a continuous
function (see \cite{JaffeTaubes}), we have near $p_j$ that
$\text{e}^{u_{\delta}}=|x-p_{j}|^{2m_{j}}f_{j,\delta}(x)$, with 
$f_{j,\delta}(x)$ a
continuous strictly positive function. Hence, (i) is established.
\par
(ii) Let $K$ be  a compact subset of
$\R^2\setminus\cup_{j\in\N}\{p_j\}$. 
In view of Lemma~\ref{lem:Ur} and 
Proposition~\ref{prop:z}, we have as $\delta \to 0^{+}$
\begin{equation}
\begin{split}
&\sup_{x \in K \cap B_{j}} 
1-\text{e}^{\varphi_j(x)U_{m_{j}}((x -p_{j})/\delta)}\le C\e^{-c/\delta}\\
&\|z_\delta(\frac{\, \cdot \,}{\delta})\|_\infty\le
C\|z_\delta\|_{\Xd}\le C\Rd
\le C\e^{-c/\delta}. 
\end{split}
\end{equation}
Therefore,  we have that for any compact
set $K \subset\R^2\setminus\cup_{j\in\N}\{p_j\}$
\begin{equation}
\label{eq:eps-del}
0\leq\sup_{x\in K} (1-\text{e}^{\ud} )\leq C \sup_{j\in\N}
\sup_{x\in K\cap B_{j}}(1-\text{e}^{\ud})\le C\text{e}^{-c/\delta}.
\end{equation}
\par
(iii) 
Let $\varphi\in C_c^\infty(\R^2)$.
Then,
\[
-\int_{\R^2}\ud\Delta\varphi=\delta^{-2}\int_{\R^2}(1-\e^u)\varphi-4\pi
m_j\varphi(p_j).
\]
We claim that
\begin{equation}
\label{eq:inttozero}
\int_{\R^2}\ud\Delta\varphi\to0
\qquad\mathrm{as\ }\delta\to0.
\end{equation}
Indeed, let $\mathrm{supp}\, \varphi\subset\cup_{k=1}^N B_{j_k}\cup
K$,
with $K$ a compact subset of $\R^2\setminus\cup_{j\in\N}\{p_j\}$.
Since $\sup_K|\ud|\le C\e^{-c/\delta}$, we have
\[
\left|\int_K\ud\Delta\varphi\right|\le
C\|\Delta\varphi\|_\infty\e^{-c/\delta}\to0.
\]
On the other hand, in view of \eqref{eq:deltasuperposition},
in $B_{j_k}$ we have $\ud(x)=U_{m_{j_k}}(|x-p_{j_k}|/\delta)+ O
(\e^{-c/\delta})$.
Note that $U_{m_{j_k}}\in L^1(\R^2)$ in view of \eqref{eq:L1}
and Lemma~\ref{lem:Ur}.
Therefore,
\begin{align*}
\sup_{1\le k\le N}\left|\int_{B_{j_k}}\ud\Delta\varphi\right|
\leq \sup_{1\le k\le N}
&\left|\int_{B_{j_k}}U_{m_{j_k}}\left(\frac{x-p_{j_k}}{\delta}\right)\Delta\varphi\right|
+ O(\e^{-c/\delta})\\
\le&\delta^2\sup_{1\le k\le N} \|\Delta\varphi\|_\infty\|U_{m_{j_k}}\|_{L^1}+
O(\e^{-c/\delta})
\le C\delta^2\to0.
\end{align*}
Hence \eqref{eq:inttozero} follows, and (iii) is established.
\end{proof}
\begin{proof}[Proof of Theorem~\ref{thm:infinite}]
For every $\delta\in(0,\delta_1)$, where $\delta_1$ 
is defined Proposition~\ref{prop:z}, we obtain a solution
$\ud$ to \eqref{eq:infinite}. 
If the $p_j$'s are doubly periodically arranged, then $\ud$
is doubly periodic in view of Lemma~\ref{lem:period}.
Furthermore, $\ud$ satisfies \eqref{eq:superposition}
in view of Lemma~\ref{lem:superposition} and of the definition
of $\delta$. Finally, $\ud$ satisfies the asymptotic behavior
as in (i)--(ii)--(iii) in view of Lemma~\ref{lem:asymp}.
Hence, Theorem~\ref{thm:infinite} is completely established.
\end{proof}

\section{Appendix}
For the reader's convenience, we sketch in this appendix the proof of
some results for
smooth, finite action critical points for the action of $\Ed$, which are relevant to our discussion.
The following results are due to Taubes~\cite{Taubes}.
Throughout this appendix all citations are referred to
the monograph of Jaffe and Taubes \cite{JaffeTaubes}.
\subsection{Derivation of equation~\eqref{eq:Taubes}}
Following \cite{JaffeTaubes} p.~53, we consider the change of variables
$A(x)=\delta^{-1}\Ap(x/\delta)$,
$\phi(x)=\phip(x/\delta)$
$\xp=x/\delta$.
We denote by $D'$, $F_{12}'$ the covariant derivative of $\Ap$ and
the curvature of $\Ap$, respectively.
Then, $\mathrm{d}A(x)=\delta^{-2}\mathrm{d}\Ap(\xp)$,
$D\phi(x)=\delta^{-1}D'\phip(\xp)$, and therefore:
\[
\int_{\R^2}\Ed(A,\phi)\,\mathrm{d}x
=\int_{\R^2}\mathcal E_1(\Ap,\phip)\,\mathrm{d}\xp,
\]
where $\mathcal E_1$ denotes $\Ed$ with $\delta=1$.
In view of Bogomol'nyi's reduction (see formula (III.1.5)),
we may rewrite the action in the form:
\begin{align*}
\int_{\R^2}\mathcal E_1&(\Ap,\phip)\,\mathrm{d}\xp\\
=&\int_{\R^2}\left\{|(D_{1}'\pm iD_{2}')\phip|^2
+(F_{12}'\pm\frac{1}{2}(|\phip|^2-1))^2\pm F_{12}'\right\}\,\mathrm{d}\xp.
\end{align*}
It follows that
\begin{align*}
\int_{\R^2}\Ed&(A,\phi)\,\mathrm{d}x\\
=&\int_{\R^2}\left\{|(D_{1}\pm iD_{2})\phi|^2
+(\delta F_{12}\pm\frac{1}{2\delta}(|\phi|^2-1))^2\pm F_{12}\right\}\,\mathrm{d}x.
\end{align*}
Here and in what follows, it is understood that
we either always choose upper signs,
or we always choose lower signs.
For smooth, finite action critical points, $(2\pi)^{-1}\int_{\R^2}F_{12}=N$ is an integer,
defining a topological class (Theorem~II.3.1 and Theorem~III.8.1).
Hence, the energy minimizers in a fixed topological class satisfy the following
the first order equations:
\begin{align}
\label{eq:sd1}
&(D_{1}\pm iD_{2})\phi=0\\
\label{eq:sd2}
&F_{12}=\pm\frac{1}{2\delta^2}(1-|\phi|^2).
\end{align}
In fact, there is no loss of generality in restricting to critical points
for the action in a given topological class (Theorem~III.10.1).
By complex analysis methods, one shows that
smooth solutions to \eqref{eq:sd1} vanish 
at most at isolated zeros of finite multiplicity.
Hence, differentiating \eqref{eq:sd1} we obtain
\[
-\Delta\ln|\phi|^2=\pm2F_{12}-4\pi\sum_{j=1}^sn_j\delta_{p_j},
\]
in the sense of distributions. 
Setting $u=\ln|\phi^2|$, we obtain from the above and \eqref{eq:sd2} 
that $u$ satisfies \eqref{eq:Taubes}.
\par
In order to define the decay properties of smooth, finite action critical points,
we define by
$u_0(x)=-\sum_{j=1}^s\ln(1+\mu|x-p_j|^{-2})$ the ``singular part" of $u$,
where $\mu>4N$.
Then, $u-u_0\in H^1(\R^2)$ (Theorem~III.3.2).
Conversely, 
if $u$ satisfies \eqref{eq:Taubes} and if $u-u_0\in H^1(\R^2)$, then
$(A,\phi)$ defined by
\begin{align*}
&\phi(z)=\exp\{\frac{1}{2}u(z)\pm i\sum_{j=1}^sm_j\arg(z-p_j)\}\\
&A_1\mp iA_2=-i(\partial_1\pm i\partial_2)\ln\phi
\end{align*}
is a smooth, finite action critical point. 
\subsection{Existence and uniqueness}
Without loss of generality we assume $\delta=1$.
The function $v=u-u_0$ satisfies the elliptic equation
with smooth coefficients
\[
-\Delta v=1-\e^{u_0+v}-g_0,
\]
where $g_0=4\sum_{j=1}^s\mu(|x-p_j|^2+\mu)^{-2}$.
Solutions in $H^1(\R^2)$ to the equation above correspond to critical points
for the functional
\[
a(v)=\int_{\R^2}\left\{\frac{1}{2}|\nabla v|^2
+(g_0-1)v+\e^{u_0}(\e^v-1)\right\},
\]
which is well-defined and differentiable on $H^1(\R^2)$.
Furthermore, $a$ is coercive and strictly convex, and therefore 
it admits a unique critical point, corresponding to the absolute minimum
(Theorem~III.4.3).
In particular, the solution to \eqref{eq:Taubes}
satisfying $u-u_0\in H^1(\R^2)$ is unique.
Finally, the critical point $(A,\phi)$ obtained from $u=u_0+v$ satisfies
the following decay estimate holds, for any $\eps>0$:
\[
|D\phi|\le\frac{3}{2}(1-|\phi|^2)\le C_\eps\e^{-(1-\eps)|x|},
\]
where $C_\eps>0$ depends on $\eps$ (Theorem~III.8.1).
 

\end{document}